# ON THE UNIFORM HYPERBOLICITY OF SOME NONUNIFORMLY HYPERBOLIC SYSTEMS

JOSÉ F. ALVES, VÍTOR ARAÚJO, AND BENOÎT SAUSSOL

ABSTRACT. We give sufficient conditions for the uniform hyperbolicity of certain nonuniformly hyperbolic dynamical systems. In particular, we show that local diffeomorphisms that are nonuniformly expanding on sets of total probability are necessarily uniformly expanding. We also present a version of this result for diffeomorphisms with nonuniformly hyperbolic sets.

## CONTENTS

## 1. INTRODUCTION

Let $f : M \to M$ be a $C^1$ local diffeomorphism of a manifold $M$ endowed with a Riemannian metric which induces a norm $|\cdot|$ on the tangent space and a volume form that we call Lebesgue measure.

The map $f$ is said to be *uniformly expanding* if there exist constants $C > 0$ and $\sigma > 1$ such that

$$(1) \qquad |df_x^n(v)| \geq C\sigma^n |v| \quad \text{for all } x \in M, \, v \in T_xM \text{ and } n \geq 1.$$

A rather complete description of the statistical properties of this kind of systems can be found in [7, 8, 11]. Recently several results have been obtained on the statistical properties of maps satisfying some weaker condition of expansion; see [3, 1]. We say that $f$ is *nonuniformly expanding at* $x \in M$ if condition (1) holds in average, meaning that

$$(2) \qquad \liminf_{n \to +\infty} \frac{1}{n} \sum_{j=0}^{n-1} \log |(df_{f^j x})^{-1}| < 0.$$

1991 *Mathematics Subject Classification.* 58F15, 58F99.

The authors were partially supported by PRODYN.

It is easy to see that if $f$ satisfies condition (1) then it also satisfies condition (2) — in fact for every $x \in M$.

From [9] one knows that $C^2$ local diffeomorphisms defined in one-dimensional manifolds are necessarily uniformly expanding if they are nonuniformly expanding almost everywhere, in the sense of Lebesgue measure. The situation is completely different in higher dimensions or with less regularity. Examples of maps which are nonuniformly expanding at Lebesgue almost every point but are not uniformly expanding are exhibited in [3].

As a by-product of our results we see that if condition (2) holds for every $x \in M$, then $f$ has to be uniformly expanding. Actually, our first theorem gives a stronger result. Before we state it, let us recall that a Borel set in a topological space is said to have *total probability* if it has probability one for every $f$ invariant probability measure.

**Theorem A.** *Let $f : M \to M$ be a $C^1$ local diffeomorphism defined in a compact Riemannian manifold. If $f$ is nonuniformly expanding on a set of total probability, then $f$ is uniformly expanding.*

We emphasize that there are examples of unimodal maps satisfying all the conditions of the theorem except the local diffeomorphism assumption, for which the conclusion fails to hold; see [5, 6].

**Corollary B.** *Let $f$ be a $C^1$ local diffeomorphism of the circle. If $f$ has a positive Lyapunov exponent on a set of total probability, then $f$ is uniformly expanding.*

We have a version of these results for $C^1$ diffeomorphisms $f \colon M \to M$ having invariant sets with some nonuniformly hyperbolic structure. Let $\Lambda \subset M$ be an $f$ invariant set with a $df$ invariant continuous splitting of the tangent bundle over $\Lambda$

$$T_\Lambda M = E^{cs} \oplus E^{cu}.$$

We say that $f$ is *nonuniformly expanding along the $E^{cu}$ direction* if, for Lebesgue almost all $x \in \Lambda$,

$$\liminf_{n \to +\infty} \frac{1}{n} \sum_{j=0}^{n-1} \log |df^{-1}|E^{cu}_{f^j x}| < 0.$$

Similarly, we say that $f$ is *nonuniformly contracting along the $E^{cs}$ direction* if, for Lebesgue almost all $x \in \Lambda$,

$$\liminf_{n \to +\infty} \frac{1}{n} \sum_{j=0}^{n-1} \log |df|E^{cs}_{f^j x}| < 0.$$

A special case is when $\Lambda$ is a *hyperbolic set*, meaning that $f$ has both uniform expansion in the $E^{cu}$ direction and uniform contraction in the $E^{cs}$ direction: there exist constants $C > 0$ and $\sigma > 1$ such that

(3) $$|df^n_x(v^u)| \geq C\sigma^n |v^u| \quad \text{and} \quad |df^n_x(v^s)| \leq C\sigma^{-n}|v^s|$$

for all $x \in M$, $v^s \in E^{cs}$, $v^u \in E^{cu}$ and $n \geq 1$.

Examples of nonuniformly hyperbolic systems which are not uniformly hyperbolic can be found in [3, 4]. The works [3, 4, 13, 2] give substantial information about the statistical properties of systems with nonuniformly hyperbolic sets.

**Theorem C.** *Let $f : M \to M$ be a $C^1$ diffeomorphism with a positively invariant set $\Lambda$ for which the tangent bundle has a continuous invariant splitting $T_\Lambda M = E^{cs} \oplus E^{cu}$. If $f$ is nonuniformly expanding (resp. contracting) along the $E^{cu}$ (resp. $E^{cs}$) direction in a set of total probability in $\Lambda$, then $f$ is uniformly expanding (resp. contracting) along the $E^{cu}$ (resp. $E^{cs}$) direction.*

In particular, if $f$ has simultaneously nonuniform expansion along the $E^{cu}$ direction and nonuniform contraction along the $E^{cs}$ direction in a set of total probability in $\Lambda$, then $\Lambda$ is a hyperbolic set.

**Corollary D.** *Let $f$ be a $C^1$ surface diffeomorphism with a positively invariant set $\Lambda$ for which the tangent bundle has a continuous invariant splitting $T_\Lambda M = E^{cs} \oplus E^{cu}$. If $f$ has a positive Lyapunov exponent in the $E^{cu}$ direction and a negative Lyapunov exponent in the $E^{cs}$ direction on a set of total probability, then $\Lambda$ is a hyperbolic set.*

Another version of these results can be given for *partially hyperbolic sets* with decomposition
$$T_\Lambda M = E^s \oplus E^c \oplus E^u,$$
where $E^s$ and $E^u$ have uniform contracting/expanding behavior.

**Corollary E.** *Let $f : M \to M$ be a $C^1$ diffeomorphism on a compact Riemannian manifold with a positively invariant set $\Lambda$ for which the tangent bundle has a continuous invariant splitting $T_\Lambda M = E^s \oplus E^c \oplus E^u$. If $E^c$ has dimension 1 and $f$ has a positive Lyapunov exponent in the $E^c$ direction or a negative Lyapunov exponent in the $E^c$ direction on a set of total probability, then $\Lambda$ is a hyperbolic set.*

**Acknowledgements:** We thank the participants in the workshop *Concepts and Techniques in Smooth Ergodic Theory* held at Imperial College, London, in July, 2001. Special thanks to the organizer, S. Luzzato, for his sponsoring of the discussions which prompted the ideias leading to this paper.

## 2. Preliminary lemmas

In this section we prove some general results for continuous dynamical systems on compact metric spaces. Let $X$ be a compact metric space, $f \colon X \to X$ a continuous function, and $\varphi \colon X \to \mathbb{R}$ also continuous such that

$$\liminf_{n \to +\infty} \frac{1}{n} \sum_{j=0}^{n-1} \varphi(f^j x) < 0 \tag{4}$$

holds in a set of points $x \in X$ with total probability.

**Lemma 2.1.** *If (4) holds in a set with total probability, then it holds for all $x \in X$.*

*Proof.* Arguing by contradiction, let us suppose that there is some $x \in X$ such that

$$\liminf_{n \to +\infty} \frac{1}{n} \sum_{j=0}^{n-1} \varphi(f^j x) \geq 0.$$

Then for every $k \in \mathbb{N}$ there must be an integer $n_k$ for which

$$\frac{1}{n_k} \sum_{j=0}^{n_k-1} \varphi(f^j x) > -\frac{1}{k}.$$

It is no restriction to assume that $n_1 < n_2 < \cdots$ and we do it. We define the sequence of probability measures

$$\mu_k = \frac{1}{n_k} \sum_{j=0}^{n_k-1} \delta_{f^j x}, \quad k \geq 1,$$

where each $\delta_{f^j x}$ is the Dirac measure on $f^j x$. Let $\mu$ be a weak* accumulation point of this sequence when $k \to +\infty$. Taking a subsequence, if necessary, we assume that $\mu_{n_k}$ converges to $\mu$. Standard arguments show that $\mu$ is $f$ invariant. Since the function $\varphi$ is continuous we have

$$\int \varphi \, d\mu = \lim_{k \to +\infty} \frac{1}{n_k} \sum_{j=0}^{n_k-1} \varphi(f^j x) \geq 0$$

by definition of $\mu$ and the way we have chosen the sequence $(n_k)_k$.

However, since we are assuming that (4) holds in a set of total probability measure, we have that $\mu$ almost every $y$ is such that

$$\tilde{\varphi}(y) = \lim_{n \to +\infty} \frac{1}{n} \sum_{j=0}^{n-1} \varphi(f^j y) < 0.$$

On the other hand,

$$\int \tilde{\varphi} \, d\mu = \int \varphi \, d\mu \geq 0$$

by the Ergodic Theorem. This gives a contradiction, thus proving the lemma. □

According to the statement of the previous lemma for each $x \in X$ there are $N(x) \in \mathbb{N}$ and $c(x) > 0$ such that

$$\frac{1}{N(x)} \sum_{j=0}^{N(x)-1} \varphi(f^j x) < -c(x).$$

Thus, by continuity, for each $x \in X$ there is a neighborhood $V_x$ of $x$ such that for every $y \in V_x$ one has

$$\frac{1}{N(x)} \sum_{j=0}^{N(x)-1} \varphi(f^j y) < -\frac{c(x)}{2}.$$

Since $X$ is compact, there is a finite cover $V_{x_1}, \ldots, V_{x_p}$ of $X$ by neighborhoods of this type. Let

(5) $$\overline{N} = \max\{N(x_1), \ldots, N(x_p)\} \quad \text{and} \quad c = \min\{c(x_1), \ldots, c(x_p)\}.$$

We define for $x \in X$
$$N_1(x) = \min\{N(x_i): x \in V_{x_i}, i = 1, \ldots, p\}$$

and a sequence of maps $N_k : M \to \{1, \ldots, \overline{N}\}$, $k \geq 0$, in the following way:

(6) $$N_0(x) = 0, \quad N_{k+1}(x) = N_k(x) + N_1(f^{N_k(x)}x).$$

(See that there is no conflict in the definition of $N_1$). For the proof of the next lemma it will be useful to take $\alpha = \max_{x \in M} \varphi(x)$

**Lemma 2.2.** *There are $C_0 > 0$ such that*
$$\sum_{j=0}^{m\overline{N}-1} \varphi(f^j x) \leq -\frac{c}{2}m + C_0$$

*for all $x \in X$ and $m \geq 1$.*

*Proof.* Let us fix $x \in X$. Given $m \geq 1$ we define
$$h = \max\{k \geq 1 : N_k(x) \leq m\overline{N}\}.$$

We must have $m\overline{N} - N_h(x) \leq \overline{N}$. It follows from (5) and (6) that for $k \geq 0$

$$\sum_{j=N_k(x)}^{N_{k+1}(x)-1} \varphi(f^j x) = \sum_{j=0}^{N_1(f^{N_k(x)})-1} \varphi(f^{N_k(x)+j}x)$$
$$\leq -\frac{c}{2}N_1(f^{N_k(x)}x)$$
$$\leq -\frac{c}{2}(N_{k+1}(x) - N_k(x))$$

Hence

$$\sum_{j=0}^{m\overline{N}-1} \varphi(f^j x) = \sum_{j=N_0(x)}^{N_1(x)-1} \varphi(f^j x) + \cdots + \sum_{j=N_{h-1}(x)}^{N_h(x)-1} \varphi(f^j x) + \sum_{j=N_h(x)}^{m\overline{N}-1} \varphi(f^j x)$$
$$\leq -\frac{c}{2}\big[(N_1(x) - N_0(x)) + \cdots + (N_h(x) - N_{h-1}(x))\big] + \alpha(m\overline{N} - N_h(x))$$
$$\leq -\frac{c}{2}N_h(x) + \alpha\overline{N}$$
$$\leq -\frac{c}{2}m + \alpha\overline{N}$$

(see that $N_h(x) \geq m\overline{N}/\overline{N} = m$). Thus we just have to choose $C_0 = \alpha\overline{N}$. $\square$

## 3. Proofs of main results

Here we use the results of the previous section to prove Theorems A and C. We start by assuming that $f: M \to M$ is a $C^1$ local diffeomorphism nonuniformly expanding on a subset of total probability. The fact that $f$ is a local diffeomorphism implies that the map

$$\lambda(x) = \log|(df_x)^{-1}|$$

is continuous on $M$. Moreover, the nonuniform expansion of $f$ implies that $\lambda$ is in the conditions of the previous section, namely it satisfies condition (4). Let $\overline{N}$, $c$ and $(N_k)_k$ be defined as in (5) and (6) for this map $\lambda$. Take also $\alpha = \max_{x \in X} \lambda(x)$.

**Proposition 3.1.** *There are constants $C > 0$ and $\sigma > 1$ such that*

$$|df_x^n(v)| \geq C\sigma^n|v|$$

*for all $x \in M$, $v \in T_xM$ and $n \geq 1$.*

*Proof.* Fix $x \in M$ and $0 \neq v \in T_xM$. Observe that

$$|(df_x^{m\overline{N}})^{-1}| \leq \exp\bigg(\sum_{j=0}^{m\overline{N}-1} \lambda(f^j x)\bigg).$$

Let $C_0 > 0$ be as in Lemma 2.2. Taking $K_0 = e^{C_0}$ and $\rho = e^{-c/2}$, for all $m \geq 1$ we have

$$|v| = |(df_x^{m\overline{N}})^{-1}(df_x^{m\overline{N}}) \cdot v| \leq K\rho^m|(df_x^{m\overline{N}}) \cdot v|,$$

which is equivalent to

$$|(df_x^{m\overline{N}}) \cdot v| \geq \frac{1}{K}(\rho^{-1})^m|v|.$$

On the other hand, we have for all $x \in M$, $v \in T_xM$ and $n \geq 1$

$$|(df_x^n) \cdot v| \geq e^{-\alpha n}|v|.$$

Let $n \geq 1$ be given. There are $m \geq 1$ and $r \in \{0, \ldots, \overline{N} - 1\}$ such that $n = m\overline{N} + r$, and so

$$\begin{aligned}
|(df_x^n) \cdot v| &= |(df_{f^{m\overline{N}}x}^r) \cdot (df_x^{m\overline{N}}) \cdot v| \\
&\geq e^{-r\alpha} \cdot |(df_x^{m\overline{N}}) \cdot v| \\
&\geq \frac{e^{-\overline{N}\alpha}}{K}(\rho^{-1})^m|v| \\
&= \frac{e^{-\overline{N}\alpha}\rho^{r/\overline{N}}}{K}(\rho^{-1/\overline{N}})^{m\overline{N}+r}|v| \\
&\geq C\sigma^n|v|.
\end{aligned}$$

It is enough to take $C = e^{-\overline{N}\alpha}\rho/K > 0$ and $\sigma = \rho^{-1/\overline{N}} > 1$. $\square$

This proposition gives precisely the conclusion of Theorem A. For proving Theorem C it is enough to take

$$\lambda^{cu}(x) = \log |df^{-1}|E_x^{cu}| \quad \text{or} \quad \lambda^{cs}(x) = \log |df|E_x^{cs}|$$

in the place of $\lambda(x)$. We use again the results of the previous section and the assumptions on $f$ to derive the desired conclusions for $df^{-1}|E_x^{cu}$ and $df|E_x^{cs}$ exactly in the same way as in Proposition 3.1.

## 4. Final remarks

We finish this paper with several remarks related to the assumptions we made for the theorems we have proved.

**Asymptotic expansion.** First we present an example which illustrates the fact that having all Lypaunov exponents positive is not enough for assuring nonuniform expansion. Consider a period 2 orbit $\{p,q\}$ for a local diffeomorphism $f$ on a surface which, for a given choice of local basis at $p$ and $q$, satisfies

$$df_p = \begin{pmatrix} 1/2 & 0 \\ 0 & 3 \end{pmatrix} \quad \text{and} \quad df_q = \begin{pmatrix} 3 & 0 \\ 0 & 1/2 \end{pmatrix}.$$

Then it is clear that both Lyapunov exponents at $p$ or $q$ are $\log(3/2) > 0$ and the limit in (2) with $x = p$ or $q$ equals $\log 2 > 0$.

**Eventual expansion.** Secondly we observe that the proof of Theorem A can be carried out under the assumption that *f is eventually expanding* at every point, meaning that for each $x \in M$ there is $N(x) \in \mathbb{N}$ such that

$$|(df_x^{N(x)})^{-1}| < 1.$$

This is essentially what we use in Lemma 2.2 and Section 3. A similar observation holds for Theorem C.

**Continuous splitting.** The final comment concerns the continuity of the splitting in Theorem C and Corollaries D and E. It is known [12, Exercice 4.1] that if a $C^1$ map $f$ admits an invariant compact set $\Lambda$ with an invariant splitting $T_\Lambda M = E^{cs} \oplus E^{cu}$ satisfying the hyperbolicity conditions (3), then this splitting is continuous. Therefore the continuity assumption on the splitting for our theorems seems natural.

Centro de Matemática da Universidade do Porto, 4099-002 Porto, Portugal
*E-mail address*: jfalves@fc.up.pt
*URL*: http://www.fc.up.pt/cmup/home/jfalves

Centro de Matemática da Universidade do Porto, 4099-002 Porto, Portugal
*E-mail address*: vdaraujo@fc.up.pt
*URL*: http://www.fc.up.pt/cmup/home/vdaraujo

LAMFA - Université de Picardie Jules Verne, 33 rue St Leu, 80039 Amiens, France
*E-mail address*: benoit.saussol@u-picardie.fr
*URL*: http://www.mathinfo.u-picardie.fr/saussol